\documentclass{amsart}

\usepackage{pdfpages}
\usepackage[procnames]{listings}
\usepackage{color}
\usepackage{cite}
\usepackage{relsize}
\usepackage{amssymb}
\numberwithin{equation}{section}
\usepackage{amsmath}
\usepackage{mathtools}
\usepackage{amssymb}
\title{The Cram\'er conjecture holds with a positive probability}
\author{Naser T. Sardari}
\date{\today}

	\newtheorem{thm}{Theorem}[section]

            \newtheorem{cram}[thm]{Cram\'er's model}

	\newtheorem{cor}[thm]{Corollary}
	
	\theoremstyle{defi}

	\theoremstyle{pf}

	\numberwithin{equation}{section}

\begin {document}
\maketitle
\begin{abstract}
We prove that a positive proportion of the intervals of any fixed scalar multiple of $\log(X)$ in the dyadic interval $[X,2X]$ contain a prime number. We also show that a positive proportion of the congruence classes modulo $q$ contain a prime number smaller than any fixed scalar multiple of  $\varphi(q)\log(q).$   
%
\end{abstract}
\tableofcontents 
\section{Introduction}
\subsection{Motivation}
Let $\pi(X)$ be the number of the prime numbers less than $X$. Then the prime number theorem states that  $\pi(X)$ is asymptotically $Li(X)=\int_{0}^{X} \frac{dt}{\log t} dt \approx \frac{X}{\log(X)}.$
Hence, on average we have one prime number in an interval of size $\log(X)$ inside the dyadic interval  $[X,2X].$ In this paper, we study how the prime numbers are distributed in 
the short intervals with the length $\lambda \log(X)$ where $\lambda>0.$ First, we give the conjectural answer that is predicted by the Cram\'er model and then the  conditional result of Gallagher proving this result by assuming the  prime $k$-tuple conjecture of Hardy and Littlewood \cite{HL}. We refer the reader to the nice exposition of Soundararajan~\cite{Soundist} for further discussion of this problem and the related results. We cite the following formulation of the  Cram\'er model from~\cite{Soundist}. 
\begin{cram}
The primes behave like independent random variables $X(n)$ $ (n \geq 3)$ with $X(n) = 1$ (the number $n$ is `prime') with probability $1/\log n$, and $X(n) = 0$ (the number $n$ is `composite') with probability $1-1/\log n.$
\end{cram}
Define $P_{k}(\lambda,X)$ to be 
\begin{equation}
P_{k}(\lambda,X):=\frac{1}{X}\# \{ X\leq n\leq 2X: \pi(n+\lambda \log(X))-\pi(n)=k \}.
\end{equation}
It follows form the above Cram\'er model that; see~\cite{Soundist} 
$$
\lim_{X\to \infty }P_{k}(\lambda,X)= \frac{e^{-\lambda}\lambda^{k} }{k!}.
$$
In fact Gallagher~  \cite[Theorem 1]{Gallagher} proved the above limit holds by assuming the  prime $k$-tuple conjecture of Hardy and Littlewood. Moreover, by using an upper bound sieve for the $k$-tuple problem Gallagher~  \cite[Theorem 2]{Gallagher}, gives an unconditional exponentially decaying upper bound for $P_{k}(\lambda,X)$ in terms of $k$.  In this paper, we   we prove the following result  

\begin{thm}\label{cramer}
Let $\lambda>0$ be any positive number and  $P_{k}(\lambda,X)$ be as above. Then, 
$$
\liminf_{X\to \infty} \sum_{k \geq 1}  P_{k}(\lambda,X)> \frac{\lambda}{4\lambda+1}+o(1)> 0.
$$
In other words, with probability at least $ \frac{\lambda}{4\lambda+1}+o(1)>0$, an interval of size $\lambda \log X$ inside the dyadic interval $[X,2X]$ contains a prime number. 
\\
\end{thm}

\begin{cor}
Cram\'er conjecture holds with a positive probability. 
\\
\end{cor}


 We also show that the $p$-adic analogue of the above statement holds that is the Linnik's conjecture. More precisely, let $q$ be an integer. We show that a positive proportion of congruence classes modulo $q$ contain a prime number smaller than $\lambda \log(q)\varphi(q)$ for any $0<\lambda$ without any assumption.

\begin{thm}\label{cramerpositive}
 Let $q$ be an integer and $X\geq \lambda \varphi(q)\log(q) $ for some fixed $\lambda >0$. Then a positive proportion, that only depends on $\lambda$, of congruence classes modulo $q$ contains a prime number smaller than   $X$. Conversely if a positive number of congruence classes modulo $q$ contains a prime number less than $X$ then $X\gg \varphi(q)\log(q) .$
 \\
\end{thm}

\begin{cor}
Linnik's conjecture holds with a positive probability. 
\end{cor}

In \cite{least}, we generalize our method to the class of the binary quadratic forms of discriminant $-D$.
We show that a positive proportion (independent of $D$) of quadratic form of discriminant $-D$ represent a prime number less than any fixed scalar multiple of $h(D)\log(D)$ by assuming a Littlewood type zero free region for the Dirichlet $L$-function $L(s,\chi_{-D})$ which holds for almost all  $D$.  

\subsection{Further questions}
An interesting question is to give unconditional  lower bounds for individual $k$
$$
\lim_{X\to \infty }P_{k}(\lambda,X)> \psi(\lambda,k)> 0,
$$
and study the limiting behavior of such lower bounds as $\lambda$ and $k$ varies.

\subsection{Acknowledgements}\noindent
I would like to thank Professor Roger Heath-Brown my mentor at MSRI for several insightful and inspiring conversations during the Spring~2017 Analytic Number Theory program at MSRI.  This material is based upon work supported by the National Science Foundation under Grant No. DMS-1440140 while the author  was in residence at the Mathematical Sciences Research Institute in Berkeley, California, during the Spring 2017 semester.

\section{Proof of Theorem~\ref{cramer}}
%
%
\begin{proof}
Let 
$$
I_n:= \pi(n+\lambda \log(X))-\pi(n),
$$
and 
$$
R_X:=\{X\leq n \leq 2X: I_n(\lambda,X)\geq 1   \}.
$$
By Cauchy inequality, we have 
\begin{equation}\label{ccramer}
R_X\big(\sum_{n=X}^{2X} I_n^2   \big)\geq \big(\sum_{n=X}^{2X} I_n \big)^2.
\end{equation}
First, we give a lower bound on the right hand side of the above inequality. By a simple double counting formula, we have 
\begin{equation}\label{moment1}
\begin{split}
\sum_{n=X}^{2X} I_n &= \lfloor \lambda \log(X)\rfloor (\pi(2X)-\pi(X))+O(\log(X)^2)
\\
& \geq (\lambda+o(1)) X.
\end{split}
\end{equation}
 Next, we give a double count formula for $S:=\sum_{n=X}^{2X} I_n^2 $. Let $H:=\{h\in\mathbb{Z}:  0 \leq h \leq  \lambda \log(X)   \}.$ $S$ gives the number of  the pairs of prime numbers $(p_1,p_2)$ such that 
 \begin{equation*}
 \begin{split}
 p_1= n+h_1
 \\
 p_2= n+h_2
 \end{split}
 \end{equation*}
 where $X<p_1,p_2<2X$ and $h_1,h_2\in H.$ We fix $h_1$ and $h_2$ and let 
 \begin{equation*}
 A(X,h_1,h_2):=\sum_{n=X}^{2X} \chi_{\mathbb{P}}(n+h_1) \chi_{\mathbb{P}}(n+h_2),
 \end{equation*}
 where $ \chi_{\mathbb{P}}$ is the characteristic function of the prime numbers. Therefore, we have 
 \begin{equation*}
 \sum_{n=X}^{2X} I_n^2 =  \sum_{\{h_1,h_2\}\subset H} A(X,h_1,h_2).
 \end{equation*}
If $h_1=h_2 $, then $A(X,h_1,h_2)=\pi(2X)-\pi(X) +O(\log(X)) $ and  the contribution of the diagonal terms are
\begin{equation}\label{diagonal}
 \sum_{h\in H} A(X,h,h)= \lfloor \lambda \log(X)\rfloor (\pi(2X)-\pi(X))+O(\log(X)^2).
\end{equation}
Next, we give an upper bound on the non-diagonal terms by applying the Selberg upper bound sieve and the result of Gallagher on the asymptotic of the average  of the Hardy-Littlewood singular series. By the Selberg upper bound sieve; see \cite[Theorem~5.7]{Halberstam}, we have 
 $$A(X,h_1,h_2)\leq 4 \mathfrak{S}(h_1,h_2)\frac{X}{\log(X)^2} \times \Big(1+O\big(\frac{\log\log(3X)+\log\log3|h_1-h_2|}{\log X}\big)   \Big),$$
 where the constant implied by $O$ term is   absolute and
\begin{equation}
\begin{split} 
\mathfrak{S}(h_1,h_2):&= \prod_{p} \big(1-\frac{\nu_{p}(h_1-h_2)}{p} \big) \big(1-\frac{1}{p}\big)^{-2},
\\
\nu_{p}(h)&=\begin{cases} 1  &\text{  if } p|h \\
2 &\text{ otherwise.}   \end{cases}
\end{split}
\end{equation}
$\mathfrak{S}(h_1,h_2)$ is the singular series in the Hardy-Littlewood prime $k$-tuples conjecture associated to the  set $\{h_1,h_2\}$.  
 By applying the above inequalities and summing over $\{h_1,h_2 \} \subset H$, we obtain
 \begin{equation}
 \begin{split}
  \sum_{n=X}^{2X} I_n^2  &=  \sum_{h_1\neq h_2 \leq \lambda \log(X)} A(X,h_1,h_2) 
  \\
  &\leq  4  \frac{X}{\log(X)^2}\sum_{h_1\neq h_2 \leq \lambda \log(X)}\mathfrak{S}(h_1,h_2).
  \end{split}
  \end{equation}
By using Gallagher's result on the average of the  Hardy-Littlewood singular series, \cite[equation (3)]{Gallagher}, we have 
$$
 \sum_{h_1\neq h_2 \leq \lambda \log(X)}\mathfrak{S}(h_1,h_2)= \big(\lambda \log(X)\big)^2(1+o(1)).
$$
 Therefore, by the above and the equation \eqref{diagonal}, we obtain 
$$
  \sum_{n=X}^{2X} I_n^2  \leq (4 \lambda^2+\lambda+ o(1)) X.
$$
By inequality \eqref{ccramer}, \eqref{moment1} and the above inequality, we obtain
$$
R_{X}\geq \big(\frac{\lambda}{4\lambda+1}+o(1) \big)X.
$$ 
This completes the proof of our theorem. 
\end{proof}

\section{Proof of Theorem~\ref{cramerpositive}}
\begin{proof}
We begin by showing that if a positive proportion of congruence classes modulo $q$ contain a prime number less than $X$, then 
$$\varphi(q)\log(q)\ll X.$$
The number of prime numbers less than $X$ is asymptotically $X/\log(X)$. We have
\begin{equation}\label{Rup}R(X,q)\leq \pi(X) \leq \frac{X}{\log(X)}.\end{equation} Assume that a positive proportion of congruence classes contain a prime number less than $X$. Then $$ c\varphi(q)<R(X,q),$$  for some positive constant $0<c$ that is independent of $q$. 
Then by the above inequality and inequality~\eqref{Rup}, we obtain
 $$c^{\prime} \varphi(q) \log(q) \leq X,$$
 for some $c^{\prime}>0$ that only depends on $c$. 
 In this section, we show that the inverse of the above necessary condition holds which means the above bound is optimal. Namely,  if $$X\geq c \varphi(q) \log(q),$$
 for some $ 0<c$ then 
 $$ R(X,q)>c^{\prime}\varphi(q),$$ for some $0<c^{\prime}\leq 1$ that only depends on $0<c$. Assume that 
\begin{equation}\label{assumption}\varphi(q)\ll X/\log(X).
\end{equation} Let $\pi(X,a,q)$ be the number of prime numbers $p<X$ such that $p\equiv a \text{ mod } q$. We proceed by applying the  Cauchy inequality and obtain 
$$R(X,q)\big(\sum_{a \text{ mod } q} \pi(X,a,q)^2\big) \geq \big(\sum_{a \text{ mod } q} \pi(X,a,q)    \big)^2.$$
 Note that $\sum_{a \text{ mod } q} \pi(X,a,q)=\pi(X)\approx X/\log(X)$, is the number of prime numbers less than $X$. Hence
 \begin{equation}\label{cauchyineq}
 R(X,q)\big(\sum_{a \text{ mod } q} \pi(X,a,q)^2\big) \geq X^2/\log(X)^2. 
 \end{equation}
 Next, we give a double count formula for $\sum_{a \text{ mod } q} \pi(X,a,q)^2$. Note that this sum counts the pairs of prime numbers $(p_1,p_2)$ such that $p_1,p_2<X$ and 
 $$p_1\equiv p_2 \text{ mod } q.$$
 This counting problem is reduced to counting the integral solutions to the following additive problem in prime numbers $(x,y)$ and some integer $t$
 \begin{equation}\label{additive}
 x=y+tq,
 \end{equation}
 such that $0<x,y<X.$ Since $0<x,y<X$, then $|tq|<X$ and hence $|t|<X/q$. Let $A(X,tq)$ denote the number of prime solutions $(x,y)$ to equation~(\ref{additive}) such that $0<x,y<X$. Then
 \begin{equation}\label{identity}
 \sum_{|t|<X/q}A(X,tq)=\sum_{a \text{ mod } q} \pi(X,a,q)^2.
 \end{equation}
  If $t=0$ then it corresponds to the diagonal elements $p_1=p_2$ and we obtain
 $$A(X,0)=\pi(X)\approx X/\log(X)$$
If $t\neq 0$ then by Selberg upper bound sieve; see \cite[Theorem~5.7]{Halberstam}, we have 
 $$A(X,tq)\leq 4 \mathfrak{S}(tq)\frac{X}{\log(X)^2} \times \Big(1+O\big(\frac{\log\log(3X)+\log\log3|tq|}{\log X}\big)   \Big),$$
 where the constant implied by $O$ term is   absolute and
\begin{equation}
\begin{split} 
\mathfrak{S}(tq):&= \prod_{p} \big(1-\frac{\nu_{p}(tq)}{p} \big) \big(1-\frac{1}{p}\big)^{-2}
\\
\nu_{p}(tq)&=\begin{cases} 1  &\text{  if } p|tq \\
2 &\text{ otherwise.}   \end{cases}
\end{split}
\end{equation}
$\mathfrak{S}(tq)$ is the singular series in the Hardy-Littlewood prime $k$-tuples conjecture associated to the  set $\{0,tq\}$.  
 By applying the above inequalities and summing over $0\leq t \leq X/q$, we obtain
 $$ \sum_{|t|<X/q}A(X,tq) \leq   X/\log(X) +4  \frac{X}{\log(X)^2}\sum_{|t|<X/q}\mathfrak{S}(tq).$$
 It follows that
 \begin{equation*}
 \begin{split}
 \mathfrak{S}(tq) &\leq \frac{q}{\varphi(q)} \mathfrak{S}(t)\prod_{p|q}(1+\frac{1}{p(p-2)}) 
 \\
 &\ll \frac{q}{\varphi(q)} \mathfrak{S}(t).
 \end{split}
 \end{equation*}
 Therefore,
 $$ \sum_{|t|<X/q}A(X,tq) \ll   X/\log(X) +  \frac{q}{\varphi(q)}\frac{X}{\log(X)^2}\sum_{|t|<X/q}\mathfrak{S}(t).$$
 By applying the result of Gallagher on the average size of the Hardy-Littlewood singular series  \cite{Gallagher}, we obtain
 \begin{equation}\label{l2}\sum_{a \text{ mod } q} \pi(X,a,q)^2 \ll X/\log(X)+ X^2/(\varphi(q)\log(X)^2).\end{equation}
 By inequalities~(\ref{l2}) and (\ref{cauchyineq}), we obtain
 $$X^2/\log(X)^2 \ll R(X,q) \big( X^2/(\varphi(q)\log(X)^2)+X/\log(X) \big).$$
 By our assumption in~(\ref{assumption}) $$\varphi(q)\leq X/\log(X).$$  Hence, $X/\log(X) \leq X^2/(\varphi(q)\log(X)^2)$ and we obtain 
 \begin{equation}
   \varphi(q) \ll R(X,d). 
 \end{equation}
 This concludes our theorem. 
 \end{proof}

\bibliographystyle{alpha}
\bibliography{cramer}

\end{document}